\documentclass{article}%
\usepackage{amsmath}
\usepackage{amsfonts}
\usepackage{amssymb}
\usepackage{graphicx}%
\providecommand{\U}[1]{\protect\rule{.1in}{.1in}}

\newenvironment{proof}[1][Proof]{\noindent\textbf{#1.} }{\ \rule{0.5em}{0.5em}}
\begin{document}

\title{On a Sturm-Liouville Type Problem with Retarded Argument}
\author{{\Large Erdo\u{g}an \c{S}en*}$^{1}${\Large , Azad Bayramov}*$^{2}$}
\date{}
\maketitle

{\Large *}$^{1}${\scriptsize {Department of Mathematics, Faculty of Arts and
Science, Nam\i k Kemal University, Tekirda\u{g}, Turkey..}}

{\scriptsize e-mail: abayramov@nku.edu.tr}$^{1},${\scriptsize esen@nku.edu.tr}%
$^{2}$



ABSTRACT. In this work a Sturm-Liouville type problem with retarded argument
which contains spectral parameter in the boundary conditions and with
transmission conditions at the point of discontinuity are investigated. We
obtained asymptotic formulas for the eigenvalues and eigenfunctions.

\section{Introduction}

We consider the boundary value problem for the differential equation
\begin{equation}
p(x)y^{\prime\prime}(x)+q(x)y(x-\Delta(x))+\lambda y(x)=0 \tag*{(1)}%
\end{equation}
on $\left[  0,\frac{\pi}{2}\right)  \cup\left(  \frac{\pi}{2},\pi\right]  ,$
with boundary conditions%

\begin{equation}
\sqrt{\lambda}y(0)+p_{1}y^{\prime}(0)=0, \tag*{(2)}%
\end{equation}%
\begin{equation}
d\lambda y(\pi)+y^{\prime}(\pi)=0, \tag*{(3)}%
\end{equation}
and transmission conditions%
\begin{equation}
\gamma_{1}y(\frac{\pi}{2}-0)=\delta_{1}y(\frac{\pi}{2}+0), \tag*{(4)}%
\end{equation}%
\begin{equation}
\gamma_{2}y^{\prime}(\frac{\pi}{2}-0)=\delta_{2}y^{\prime}(\frac{\pi}{2}+0),
\tag*{(5)}%
\end{equation}
where $p(x)=p_{1}^{2}$ if $x\in\left[  0,\frac{\pi}{2}\right)  $ and
$p(x)=p_{2}^{2}$ if $x\in\left(  \frac{\pi}{2},\pi\right]  $, the real-valued
function $q(x)$ is continuous in $\left[  0,\frac{\pi}{2}\right)  \cup\left(
\frac{\pi}{2},\pi\right]  ~$and has a finite limit $q(\frac{\pi}{2}\pm
0)=\lim_{x\rightarrow\frac{\pi}{2}\pm0}q(x),$ the real valued function
$\Delta(x)\geq0$ continuous in $\left[  0,\frac{\pi}{2}\right)  \cup\left(
\frac{\pi}{2},\pi\right]  $ and has a finite limit $\Delta(\frac{\pi}{2}%
\pm0)=\lim_{x\rightarrow\frac{\pi}{2}\pm0}\Delta(x)$, $x-\Delta(x)\geq0,$ if
$\ x\in\left[  0,\frac{\pi}{2}\right)  ;x-\Delta(x)\geq\frac{\pi}{2},$ if
$x\in\left(  \frac{\pi}{2},\pi\right]  ;$ $\lambda$ is a real spectral
parameter; $p_{1},p_{2},\gamma_{1},\gamma_{2},\delta_{1},\delta_{2},d$ are
arbitrary real numbers and $\left\vert \gamma_{i}\right\vert +\left\vert
\delta_{i}\right\vert \neq0$ for $i=1,2$. Also $\gamma_{1}\delta_{2}%
p_{1}=\gamma_{2}\delta_{1}p_{2}$ holds.

The present article is devoted to studying the properties of the eigenvalues
and eigenfunctions of the boundary value problem (1)-(5). Boundary value
problems for ordinary differential equations with a parameter in the boundary
conditions in various statements were studied in many articles $[1-10]$.

The article $[11]$ is devoted to the study of the asymptotics of the solutions
to the Sturm-Liouville problem with the potential and the spectral parameter
having discontinuity of the first kind in the domain of definition of the solution.

Boundary value problems for differential equations of the second order with
retarded argument were studied in $[12,13]$.

The asymptotic formulas for the eigenvalues and eigenfunctions of boundary
problem of discontinuous Sturm-Liouville problem with transmission conditions
and with the boundary conditions which include spectral parameter were
obtained in $\left[  14,15\right]  $.

In the present considered problem's both boundary conditions involves spectral
parameter. The main result of the present paper is Theorem 4 and Theorem 5 on
asymptotic formulas for eigenvalues and eigenfunctions.

It must be also noted that some problems with transmission conditions which
arise in mechanics (thermal condition problem for a thin laminated plate) were
studied in $\left[  16\right]  $.

Let $w_{1}(x,\lambda)$ be a solution of Eq. (1) on $\left[  0,\frac{\pi}%
{2}\right]  ,$ satisfying the initial conditions%
\begin{equation}
w_{1}\left(  0,\lambda\right)  =p_{1},w_{1}^{\prime}\left(  0,\lambda\right)
=-\sqrt{\lambda} \tag*{(6)}%
\end{equation}
The conditions $(6)$ define a unique solution of Eq. $(1)$ on $\left[
0,\frac{\pi}{2}\right]  $($\left[  13\right]  $, p. 12).

After defining above solution we shall define the solution $w_{2}\left(
x,\lambda\right)  $ of Eq. (1) on $\left[  \frac{\pi}{2},\pi\right]  $ by
means of the solution $w_{1}\left(  x,\lambda\right)  $ by the initial
conditions%
\begin{equation}
w_{2}\left(  \frac{\pi}{2},\lambda\right)  =\gamma_{1}\delta_{1}^{-1}%
w_{1}\left(  \frac{\pi}{2},\lambda\right)  ,\quad\omega_{2}^{\prime}%
({\frac{\pi}{2}},\>\lambda)=\gamma_{2}\delta_{2}^{-1}\omega_{1}^{\prime
}({\frac{\pi}{2}},\>\lambda). \tag*{(7)}%
\end{equation}
The conditions (7) are defined as a unique solution of Eq. (1) on $\left[
\frac{\pi}{2},\pi\right]  .$

Consequently, the function $w\left(  x,\lambda\right)  $ is defined on
$\left[  0,\frac{\pi}{2}\right)  \cup\left(  \frac{\pi}{2},\pi\right]  $ by
the equality%
\[
w(x,\lambda)=\left\{
\begin{array}
[c]{ll}%
\omega_{1}(x,\lambda), & x\in\lbrack0,{\frac{\pi}{2}})\\
\omega_{2}(x,\lambda), & x\in({\frac{\pi}{2}},\pi]
\end{array}
\right.
\]
is a such solution of the Eq. (1) on $\left[  0,\frac{\pi}{2}\right)
\cup\left(  \frac{\pi}{2},\pi\right]  ;$which satisfies boundary conditions
and transmission conditions.

\textbf{Lemma 1} \ Let $w\left(  x,\lambda\right)  $ be a solution of Eq.$(1)$
and\ $\lambda>0.$ Then the following integral equations hold:\
\begin{align}
w_{1}(x,\lambda)  &  =\sqrt{2}p_{1}\cos\left(  \frac{s}{p_{1}}x+\frac{\pi}%
{4}\right) \nonumber\\
&  -\frac{1}{s}\int\limits_{0}^{{x}}\frac{q\left(  \tau\right)  }{p_{1}}%
\sin\frac{s}{p_{1}}\left(  x-\tau\right)  w_{1}\left(  \tau-\Delta\left(
\tau\right)  ,\lambda\right)  d\tau\text{ \ }\left(  s=\sqrt{\lambda}%
,\lambda>0\right)  , \tag{8}%
\end{align}%
\begin{align}
w_{2}(x,\lambda)  &  =\frac{\gamma_{1}}{\delta_{1}}w_{1}\left(  \frac{\pi}%
{2},\lambda\right)  \cos\frac{s}{p_{2}}\left(  x-\frac{\pi}{2}\right)
+\frac{\gamma_{2}p_{2}w_{1}^{\prime}\left(  \frac{\pi}{2},\lambda\right)
}{s\delta_{2}}\sin\frac{s}{p_{2}}\left(  x-\frac{\pi}{2}\right) \nonumber\\
&  -\frac{1}{s}\int\limits_{\pi/2}^{{x}}\frac{q\left(  \tau\right)  }{p_{2}%
}\sin\frac{s}{p_{2}}\left(  x-\tau\right)  w_{2}\left(  \tau-\Delta\left(
\tau\right)  ,\lambda\right)  d\tau\text{ \ }\left(  s=\sqrt{\lambda}%
,\lambda>0\right)  . \tag{9}%
\end{align}

\begin{proof}
To prove this, it is enough to substitute $\>-\frac{s^{2}}{p_{1}^{2}}%
\omega_{1}(\tau,\lambda)-\omega_{1}^{\prime\prime}(\tau,\lambda)\>$ and
$\>-\frac{s^{2}}{p_{2}^{2}}\omega_{2}(\tau,\lambda)-\omega_{2}^{\prime\prime
}(\tau,\lambda)\>$ instead of $\>-\frac{q(\tau)}{p_{1}^{2}}\omega_{1}%
(\tau-\Delta(\tau),\lambda)\>$ and $\>-\frac{q(\tau)}{p_{2}^{2}}\omega
_{2}(\tau-\Delta(\tau),\lambda)\>$ in the integrals in (8) and (9)
respectively and integrate by parts twice.
\end{proof}

\textbf{Theorem 1 \ }The problem $(1)-(5)$ can have only simple eigenvalues.

\begin{proof}
Let $\widetilde{\lambda}$ be an eigenvalue of the problem $(1)-(5)$ and%
\[
\widetilde{u}(x,\widetilde{\lambda})=\left\{
\begin{array}
[c]{ll}%
\widetilde{u}_{1}(x,\widetilde{\lambda}), & x\in\lbrack0,{\frac{\pi}{2}}),\\
\widetilde{u}_{2}(x,\widetilde{\lambda}), & x\in({\frac{\pi}{2}},\pi]
\end{array}
\right.
\]
be a corresponding eigenfunction. Then from $(2)$ and $(6)$ it follows that
the determinant%
\[
W\left[  \widetilde{u}_{1}(0,\widetilde{\lambda}),w_{1}(0,\widetilde{\lambda
})\right]  =\left\vert
\begin{array}
[c]{c}%
\widetilde{u}_{1}(0,\widetilde{\lambda})\text{ \ \ \ \ \ \ }p_{1}\\
\widetilde{u}_{1}^{\prime}(0,\widetilde{\lambda})\text{ \ \ }-\sqrt{\lambda}%
\end{array}
\right\vert =0,
\]
and by Theorem 2.2.2 in $\left[  13\right]  $ the functions $\widetilde{u}%
_{1}(x,\widetilde{\lambda})$ and $w_{1}(x,\widetilde{\lambda})$ are linearly
dependent on $\left[  0,\frac{\pi}{2}\right]  $. We can also prove that the
functions $\widetilde{u}_{2}(x,\widetilde{\lambda})$ and $w_{2}(x,\widetilde
{\lambda})$ are linearly dependent on $\left[  \frac{\pi}{2},\pi\right]  $.
Hence%
\begin{equation}
\widetilde{u}_{1}(x,\widetilde{\lambda})=K_{i}w_{i}(x,\widetilde{\lambda
})\text{ \ \ \ }\left(  i=1,2\right)  \tag*{(10)}%
\end{equation}
for some $K_{1}\neq0$ and $K_{2}\neq0$. We must show that $K_{1}=K_{2}$.
Suppose that $K_{1}\neq K_{2}$. From the equalities $(4)$ and $(10)$, we have%
\begin{align*}
\gamma_{1}\widetilde{u}(\frac{\pi}{2}-0,\widetilde{\lambda})-\delta
_{1}\widetilde{u}(\frac{\pi}{2}+0,\widetilde{\lambda})  &  =\gamma
_{1}\widetilde{u_{1}}(\frac{\pi}{2},\widetilde{\lambda})-\delta_{1}%
\widetilde{u_{2}}(\frac{\pi}{2},\widetilde{\lambda})\\
&  =\gamma_{1}K_{1}w_{1}(\frac{\pi}{2},\widetilde{\lambda})-\delta_{1}%
K_{2}w_{2}(\frac{\pi}{2},\widetilde{\lambda})\\
&  =\gamma_{1}K_{1}\delta_{1}\gamma_{1}^{-1}w_{2}(\frac{\pi}{2},\widetilde
{\lambda})-\delta_{1}K_{2}w_{2}(\frac{\pi}{2},\widetilde{\lambda})\\
&  =\delta_{1}\left(  K_{1}-K_{2}\right)  w_{2}(\frac{\pi}{2},\widetilde
{\lambda})=0.
\end{align*}
Since $\delta_{1}\left(  K_{1}-K_{2}\right)  \neq0$ it follows that
\begin{equation}
w_{2}\left(  \frac{\pi}{2},\widetilde{\lambda}\right)  =0. \tag{11}%
\end{equation}
By the same procedure from equality $(5)$ we can derive that%
\begin{equation}
w_{2}^{^{\prime}}\left(  \frac{\pi}{2},\widetilde{\lambda}\right)  =0.
\tag{12}%
\end{equation}
From the fact that $w_{2}(x,\widetilde{\lambda})$ is a solution of the
differential Eq. $(1)$ on $\left[  \frac{\pi}{2},\pi\right]  $ and satisfies
the initial conditions $(11)$ and $(12)$ it follows that $w_{1}(x,\widetilde
{\lambda})=0$ identically on $\left[  \frac{\pi}{2},\pi\right]  $ (cf. [13, p.
12, Theorem 1.2.1]).

By using we may also find%
\[
w_{1}\left(  \frac{\pi}{2},\widetilde{\lambda}\right)  =w_{1}^{^{\prime}%
}\left(  \frac{\pi}{2},\widetilde{\lambda}\right)  =0.
\]
From the latter discussions of $w_{2}(x,\widetilde{\lambda})$ it follows that
$w_{1}(x,\widetilde{\lambda})=0$ identically on $\left[  0,\frac{\pi}%
{2}\right)  \cup\left(  \frac{\pi}{2},\pi\right]  $. But this contradicts
$(6)$, thus completing the proof.
\end{proof}

\section{An existance theorem}

The function $\omega(x,\>\lambda)\>$defined in section $1$ is a nontrivial
solution of Eq. $(1)$ satisfying conditions $(2),(4)$ and $(5)$.
Putting$\>\omega(x,\>\lambda)\>$into $(3)$, we get the characteristic equation%
\begin{equation}
F(\lambda)\equiv\omega^{\prime}(\pi,\>\lambda)+d\lambda\omega(\pi
,\>\lambda)=0. \tag*{(13)}%
\end{equation}

By Theorem 1 the set of eigenvalues of boundary-value problem (1)-(5)
coincides with the set of real roots of Eq. (13). Let $\>q_{1}=\frac{1}{p_{1}%
}\int\limits_{0}^{{\pi/2}}|q(\tau)|d\tau$ and $q_{2}=\frac{1}{p_{2}}%
\int\limits_{{\pi/2}}^{{\pi}}q(\tau)d\tau$

\textbf{Lemma 2} \ $(1)$ Let $\lambda\geq4q_{1}^{2}$. Then for the solution
$w_{1}\left(  x,\lambda\right)  $ of Eq. $(8)$, the following inequality
holds$:$%
\begin{equation}
\left\vert w_{1}\left(  x,\lambda\right)  \right\vert \leq2\sqrt{2}\left\vert
p_{1}\right\vert ,\text{ \ \ }x\in\left[  0,\frac{\pi}{2}\right]  .
\tag*{(14)}%
\end{equation}
$(2)$ Let $\lambda\geq\max\left\{  4q_{1}^{2},4q_{2}^{2}\right\}  $. Then for
the solution $w_{2}\left(  x,\lambda\right)  $ of Eq. $(9)$, the following
inequality holds$:$%
\begin{equation}
\left\vert w_{2}\left(  x,\lambda\right)  \right\vert \leq4\sqrt{2}\left\vert
p_{1}\right\vert \left\{  \left\vert \frac{\gamma_{1}}{\delta_{1}}\right\vert
+\left\vert \frac{p_{2}\gamma_{2}}{4p_{1}\delta_{2}}\right\vert \right\}
,\text{\ \ }x\in\left[  \frac{\pi}{2},\pi\right]  \tag*{(15)}%
\end{equation}

\begin{proof}
Let $B_{1\lambda}=\max_{\left[  0,\frac{\pi}{2}\right]  }\left\vert
w_{1}\left(  x,\lambda\right)  \right\vert $. Then from $(8)$, it follows
that, for every $\lambda>0$, the following inequality holds:%
\[
B_{1\lambda}\leq\sqrt{2}p_{1}+\frac{1}{s}B_{1\lambda}q_{1}.
\]
If $s\geq2q_{1}$ we get $(14)$. Differentiating $(8)$ with respect to $x$, we
have%
\begin{equation}
w_{1}^{\prime}(x,\lambda)=-s\sqrt{2}\sin\left(  \frac{sx}{p_{1}}+\frac{\pi}%
{4}\right)  -\frac{1}{p_{1}^{2}}\int\limits_{0}^{x}q(\tau)\cos\frac{s}{p_{1}%
}\left(  x-\tau\right)  w_{1}(\tau-\Delta\left(  \tau\right)  ,\lambda
)d\tau\tag{16}%
\end{equation}
From $(16)$ and $(14)$, it follows that, for $s\geq2q_{1}$, the following
inequality holds:%
\[
\left\vert w_{1}^{\prime}(x,\lambda)\right\vert \leq s\sqrt{2}+2\sqrt
{2}\left\vert q_{1}\right\vert .
\]
Hence%
\begin{equation}
\frac{\left\vert w_{1}^{\prime}(x,\lambda)\right\vert }{s}\leq\sqrt{2}
\tag{17}%
\end{equation}
Let $B_{2\lambda}=\max_{\left[  \frac{\pi}{2},\pi\right]  }\left\vert
w_{2}\left(  x,\lambda\right)  \right\vert $. Then from $(9),(14)$ and $(17)$
it follows that, for $s\geq2q_{1}$, the following inequalities holds:%
\begin{align*}
B_{2\lambda}  &  \leq4\sqrt{2}\left\vert \frac{p_{1}\gamma_{1}}{\delta_{1}%
}\right\vert +\left\vert \frac{p_{2}\gamma_{2}}{\delta_{2}}\right\vert
\frac{1}{\left\vert p_{1}q_{1}\right\vert }\sqrt{2}+\frac{B_{2\lambda}}{2},\\
B_{2\lambda}  &  \leq4\sqrt{2}\left\vert p_{1}\right\vert \left\{  \left\vert
\frac{\gamma_{1}}{\delta_{1}}\right\vert +\left\vert \frac{p_{2}\gamma_{2}%
}{4p_{1}\delta_{2}}\right\vert \right\}  .
\end{align*}
Hence if $\lambda\geq\max\left\{  4q_{1}^{2},4q_{2}^{2}\right\}  $ we get
$(15)$.
\end{proof}

\textbf{Theorem 2} \ The problem $(1)-(5)$ has an infinite set of positive eigenvalues.

\begin{proof}
Differentiating $(9)$ with respect to$\>x$, we get%
\begin{align}
w_{2}^{\prime}(x,\lambda)  &  =-\frac{s\gamma_{1}}{p_{2}\delta_{1}}%
w_{1}\left(  \frac{\pi}{2},\lambda\right)  \sin\frac{s}{p_{2}}\left(
x-\frac{\pi}{2}\right)  +\frac{\gamma_{2}w_{1}^{\prime}\left(  \frac{\pi}%
{2},\lambda\right)  }{\delta_{2}}\cos\frac{s}{p_{2}}\left(  x-\frac{\pi}%
{2}\right) \nonumber\\
&  -\frac{1}{p_{2}^{2}}\int\limits_{\pi/2}^{{x}}q(\tau)\cos\frac{s}{p_{2}%
}\left(  x-\tau\right)  w_{2}(\tau-\Delta\left(  \tau\right)  ,\lambda)d\tau.
\tag{18}%
\end{align}
From $(8),(9),(13),(16)$ and $(18)$, we get%
\begin{equation}
-{\frac{s\gamma_{1}}{p_{2}\delta_{1}}}\biggl(\sqrt{2}p_{1}\cos\left(
\frac{s\pi}{2p_{1}}+\frac{\pi}{4}\right)  -{\frac{1}{sp_{1}}}\int
\limits_{0}^{{\frac{\pi}{2}}}q(\tau)\sin\frac{s}{p_{1}}({\frac{\pi}{2}}%
-\tau)\omega_{1}(\tau-\Delta(\tau),\lambda)d\tau\biggr)\nonumber
\end{equation}%
\[
\times\sin{\frac{s\pi}{2p_{2}}}%
\]%
\[
+{\frac{\gamma_{2}}{\delta_{2}}}\biggl(-s\sqrt{2}\sin\left(  \frac{s\pi
}{2p_{1}}+\frac{\pi}{4}\right)  -\frac{1}{p_{1}^{2}}\int\limits_{0}%
^{{\frac{\pi}{2}}}q(\tau)\cos\frac{s}{p_{1}}({\frac{\pi}{2}}-\tau)\omega
_{1}(\tau-\Delta(\tau),\lambda)d\tau\biggr)
\]%
\[
\times\cos{\frac{s\pi}{2p_{2}}}-\frac{1}{p_{2}^{2}}\int\limits_{\pi/2}^{{\pi}%
}q(\tau)\cos\frac{s}{p_{2}}({\pi}-\tau)\omega_{2}(\tau-\Delta(\tau
),\lambda)d\tau
\]%
\[
+\lambda d\left(  \frac{\gamma_{1}}{\delta_{1}}\left[  \sqrt{2}p_{1}%
\cos\left(  \frac{s\pi}{2p_{1}}+\frac{\pi}{4}\right)  {-\frac{1}{sp_{1}}}%
\int\limits_{0}^{{\frac{\pi}{2}}}q(\tau)\sin\frac{s}{p_{1}}({\frac{\pi}{2}%
}-\tau)\omega_{1}(\tau-\Delta(\tau),\lambda)d\tau\right]  \right.
\]%
\[
\times\cos{\frac{s\pi}{2p_{2}}}%
\]%
\[
+\frac{\gamma_{2}p_{2}}{\delta_{2}s}\left[  -s\sqrt{2}\sin\left(  \frac{s\pi
}{2p_{1}}+\frac{\pi}{4}\right)  {-}\frac{1}{p_{1}^{2}}\int\limits_{0}%
^{{\frac{\pi}{2}}}q(\tau)\cos\frac{s}{p_{1}}({\frac{\pi}{2}}-\tau)\omega
_{1}(\tau-\Delta(\tau),\lambda)d\tau\right]
\]%
\begin{equation}
\times\sin{\frac{s\pi}{2p_{2}}}\left.  -{\frac{1}{sp_{2}}}\int\limits_{\frac
{\pi}{2}}^{{\pi}}q(\tau)\sin\frac{s}{p_{2}}({\pi}-\tau)\omega_{2}(\tau
-\Delta(\tau),\lambda)d\tau\right)  =0 \tag{19}%
\end{equation}
Let $\>\lambda\>$ be sufficiently large. Then, by $(14)$ and $(15)$, the Eq.
$(19)$ may be rewritten in the form%
\begin{equation}
s\cos\left(  s\pi\frac{p_{1}+p_{2}}{2p_{1}p_{2}}+\frac{\pi}{4}\right)  +O(1)=0
\tag{20}%
\end{equation}
Obviously, for large$\>s\>$Eq. $(20)$ has an infinite set of roots. Thus the
theorem is proved.
\end{proof}

\section{Asymptotic Formulas for Eigenvalues and Eigenfunctions}

Now we begin to study asymptotic properties of eigenvalues and eigenfunctions.
In the following we shall assume that$\>s\>$is sufficiently large. From $(8)$
and $(14)$, we get
\[
\omega_{1}(x,\>\lambda)=O(1)\quad\mbox{on}\quad\lbrack0,\>{\frac{\pi}{2}%
}]\eqno(21)
\]
From $(9)$ and $(15)$, we get
\[
\omega_{2}(x,\>\lambda)=O(1)\quad\mbox{on}\quad\lbrack{\frac{\pi}{2}}%
,\>\pi]\eqno(22)
\]
The existence and continuity of the derivatives $\>\omega_{1s}^{\prime
}(x,\>\lambda)\>$for $\>0\leq x\leq{\frac{\pi}{2}},\>|\lambda|<\infty$,
and$\>\omega_{2s}^{\prime}(x,\>\lambda)\>$for $\>{\frac{\pi}{2}}\leq x\leq
\pi,\>|\lambda|<\infty$, follows from Theorem 1.4.1 in \cite{Nk}.
\[
\>\omega_{1s}^{\prime}(x,\>\lambda)=O(1),\quad x\in\lbrack0,\>{\frac{\pi}{2}%
}]\quad\mbox{and}\quad\omega_{2s}^{\prime}(x,\>\lambda)=O(1),\quad x\in
\lbrack\>{\frac{\pi}{2}},\>\pi].\eqno(23)
\]

\textbf{Theorem 3} \ Let $n$ be a natural number. For each sufficiently large
$n,$ in case $1,$ there is exactly one eigenvalue of the problem
$(1)-(5)$\ near $\frac{p_{1}^{2}p_{2}^{2}}{4(p_{1}+p_{2})^{2}}\left(
4n-3\right)  ^{2}$

\begin{proof}
We consider the expression which is denoted by$\>O(1)$ in the Eq. $(20)$. If
formulas $(21)-(23)$ are taken into consideration, it can be shown by
differentiation with respect to$\>s$ that for large$\>s\>$this expression has
bounded derivative.It is obvious that for large$\>s$ the roots of Eq. $(20)$
are situated close to entire numbers. We shall show that, for large$\>n$, only
one root $(20)$ lies near to each $\frac{p_{1}^{2}p_{2}^{2}}{4(p_{1}%
+p_{2})^{2}}\left(  4n-3\right)  ^{2}$. We consider the function
$\>\phi(s)=s\cos\left(  s\pi\frac{p_{1}+p_{2}}{2p_{1}p_{2}}+\frac{\pi}%
{4}\right)  +O(1)\>$. Its derivative, which has the form
\[
\phi^{\prime}(s)=\cos\left(  s\pi\frac{p_{1}+p_{2}}{2p_{1}p_{2}}+\frac{\pi}%
{4}\right)  -s\pi\frac{p_{1}+p_{2}}{2p_{1}p_{2}}\sin\left(  s\pi\frac
{p_{1}+p_{2}}{2p_{1}p_{2}}+\frac{\pi}{4}\right)  +O(1),
\]
does not vanish for$\>s\>$close to$\>n\>$for sufficiently large$\>n\>$. Thus
our assertion follows by Rolle's Theorem.
\end{proof}

Let $\>n\>$ be sufficiently large. In what follows we shall denote
by$\>\lambda_{n}=s_{n}^{2}$ the eigenvalue of the problem $(1)-(5)$ situated
near $\frac{p_{1}^{2}p_{2}^{2}}{4(p_{1}+p_{2})^{2}}\left(  4n-3\right)  ^{2}$.
We set $s_{n}=\frac{p_{1}p_{2}\left(  4n-3\right)  }{2\left(  p_{1}%
+p_{2}\right)  }+\delta_{n}$. From $(20)$ it follows that $\delta_{n}=O\left(
\frac{1}{n}\right)  $. Consequently%
\[
s_{n}=\frac{p_{1}p_{2}\left(  4n-3\right)  }{2\left(  p_{1}+p_{2}\right)
}+O\bigl ({\frac{1}{n}}\bigr ).\eqno(24)
\]
The formula $(24)$ make it possible to obtain asymptotic expressions for
eigenfunction of the problem $(1)-(5)$.
From $(8),(16)$ and $(21)$, we get
\[
\omega_{1}(x,\>\lambda)=\sqrt{2}p_{1}\cos\left(  \frac{sx}{p_{1}}+\frac{\pi
}{4}\right)  +O\bigl ({\frac{1}{s}}\bigr ),\eqno(25)
\]%
\[
\omega_{1}^{^{\prime}}(x,\>\lambda)=-\sqrt{2}s\sin\left(  \frac{sx}{p_{1}%
}+\frac{\pi}{4}\right)  +O\bigl ({1}\bigr ).\eqno(26)
\]
From $(9),(22),(25)$ and $(26)$, we get%
\[
\omega_{2}(x,\>\lambda)={\frac{\sqrt{2}\gamma_{1}{p}_{1}}{\delta_{1}}}%
\cos\left(  \frac{s\pi(p_{2-}p_{1})}{2p_{1}p_{2}}+\frac{sx}{p_{2}}+\frac{\pi
}{4}\right)  +O\bigl ({\frac{1}{s}}\bigr )\eqno(27)
\]
By putting $(24)$ in the $(25)$ and $(27)$, we derive that%
\begin{align*}
u_{1n}  &  =w_{1}\left(  x,\lambda_{n}\right)  =\sqrt{2}{p}_{1}\cos\left(
\frac{p_{2}\left(  4n-3\right)  x}{2\left(  p_{1}+p_{2}\right)  }+\frac{\pi
}{4}\right)  +O\bigl ({\frac{1}{n}}\bigr ),\\
u_{2n}  &  =w_{2}\left(  x,\lambda_{n}\right)  ={\frac{\sqrt{2}{p}_{1}{\gamma
}_{1}}{\delta_{1}}}\cos\left(  \frac{p_{1}\left(  4n-3\right)  x}{2\left(
p_{1}+p_{2}\right)  }+\frac{\pi}{4}\left(  1+\frac{\left(  p_{2}-p_{1}\right)
(4n-3)}{4\left(  p_{1}+p_{2}\right)  }\right)  \right)  +O\bigl ({\frac{1}{n}%
}\bigr ).
\end{align*}
Hence the eigenfunctions$\>u_{n}(x)\>$have the following asymptotic
representation:%
\[
u_{n}(x)=\left\{
\begin{array}
[c]{lll}%
\sqrt{2}{p}_{1}\cos\left(  \frac{p_{2}\left(  4n-3\right)  x}{2\left(
p_{1}+p_{2}\right)  }+\frac{\pi}{4}\right)  +O\bigl ({\frac{1}{n}}\bigr ) &
\mbox{for} & x\in\lbrack0,{\frac{\pi}{2}}),\\
{\frac{\sqrt{2}{p}_{1}{\gamma}_{1}}{\delta_{1}}}\cos\left(  \frac{p_{1}\left(
4n-3\right)  x}{2\left(  p_{1}+p_{2}\right)  }+\frac{\pi}{4}\left(
1+\frac{\left(  p_{2}-p_{1}\right)  (4n-3)}{4\left(  p_{1}+p_{2}\right)
}\right)  \right)  +O\bigl ({\frac{1}{n}}\bigr ) & \mbox{for} & x\in
({\frac{\pi}{2}},\pi].
\end{array}
\right.
\]
Under some additional conditions the more exact asymptotic formulas which
depend upon the retardation may be obtained. Let us assume that the following
conditions are fulfilled:

\noindent\textbf{a)} The derivatives $\>q^{\prime}(x)\>$ and $\>\Delta
^{\prime\prime}(x)\>$ exist and are bounded in$\>[0,{\frac{\pi}{2}}%
)\bigcup({\frac{\pi}{2}},\pi]\>$ and have finite limits $\>q^{\prime}%
({\frac{\pi}{2}}\pm0)=\lim\limits_{x\rightarrow{\frac{\pi}{2}}\pm0}q^{\prime
}(x)\>$ and $\>\Delta^{\prime\prime}({\frac{\pi}{2}}\pm0)=\lim
\limits_{x\rightarrow{\frac{\pi}{2}}\pm0}\Delta^{\prime\prime}(x)$, respectively.

\noindent\textbf{b)} $\>\Delta^{\prime}(x)\leq1 \>$ in$\>[0,{\frac{\pi}{2}%
})\bigcup({\frac{\pi}{2}},\pi]$, $\>\Delta(0)=0\>$ and $\>\lim
\limits_{x\rightarrow{\frac{\pi}{2}}+ 0}\Delta(x)=0$.

By using b), we have
\[
x-\Delta(x)\geq0,\quad\mbox{for}\>x\in\lbrack0,{\frac{\pi}{2}})\quad
\mbox{and}\quad x-\Delta(x)\geq{\frac{\pi}{2}},\quad\mbox{for}\>x\in
({\frac{\pi}{2}},\pi]\eqno(28)
\]

From $(25)$, $(27)$ and $(28)$ we have%
\begin{equation}
w_{1}\left(  \tau-\Delta\left(  \tau\right)  ,\lambda\right)  =\sqrt{2}%
p_{1}\cos\left(  \frac{s}{p_{1}}\left(  \tau-\Delta\left(  \tau\right)
\right)  +\frac{\pi}{4}\right)  +O\bigl ({\frac{1}{s}}\bigr ), \tag{29}%
\end{equation}%
\begin{equation}
w_{2}\left(  \tau-\Delta\left(  \tau\right)  ,\lambda\right)  ={\frac{\sqrt
{2}p_{1}{\gamma}_{1}}{\delta_{1}}}\cos\left(  \frac{s\pi(p_{2-}p_{1})}%
{2p_{1}p_{2}}+\frac{s\left(  \tau-\Delta\left(  \tau\right)  \right)  }{p_{2}%
}+\frac{\pi}{4}\right)  +O\bigl ({\frac{1}{s}}\bigr ). \tag{30}%
\end{equation}
Putting these expressions into $(19)$, we have%
\begin{align*}
0  &  =-\frac{\gamma_{2}}{\delta_{2}}\sin\left(  \frac{s\pi\left(  p_{1}%
+p_{2}\right)  }{2p_{1}p_{2}}+\frac{\pi}{4}\right)  +\frac{sd{\gamma}_{1}%
p_{1}}{\delta_{1}}\cos\left(  \frac{s\pi\left(  p_{1}+p_{2}\right)  }%
{2p_{1}p_{2}}+\frac{\pi}{4}\right) \\
&  -\frac{d{\gamma}_{1}}{p_{1}\delta_{1}}\left[  \sin\left(  \frac{s\pi\left(
p_{1}+p_{2}\right)  }{2p_{1}p_{2}}\right)  \int\limits_{0}^{\pi/2}\frac
{\sqrt{2}q\left(  \tau\right)  }{2}\cos\left(  \frac{s\left(  2\tau
-\Delta(\tau)\right)  }{p_{1}}+\frac{\pi}{4}\right)  d\tau\right. \\
&  -\cos\left(  \frac{s\pi\left(  p_{1}+p_{2}\right)  }{2p_{1}p_{2}}\right)
\int\limits_{0}^{\pi/2}\frac{\sqrt{2}q\left(  \tau\right)  }{2}\sin\left(
\frac{s\left(  2\tau-\Delta(\tau)\right)  }{p_{1}}+\frac{\pi}{4}\right)
d\tau\\
&  +\sin\left(  \frac{s\pi\left(  p_{1}+p_{2}\right)  }{2p_{1}p_{2}}\right)
\sin\int\limits_{0}^{\pi/2}\frac{\sqrt{2}q\left(  \tau\right)  }{2}\cos\left(
\frac{s\Delta(\tau)}{p_{1}}-\frac{\pi}{4}\right)  d\tau\\
&  \left.  -\cos\left(  \frac{s\pi\left(  p_{1}+p_{2}\right)  }{2p_{1}p_{2}%
}\right)  \int\limits_{0}^{\pi/2}\frac{\sqrt{2}q\left(  \tau\right)  }{2}%
\sin\left(  \frac{s\Delta(\tau)}{p_{1}}-\frac{\pi}{4}\right)  d\tau\right] \\
&  -\frac{d\sqrt{2}{\gamma}_{1}}{p_{2}\delta_{1}}\left[  \sin\left(  s\pi
\frac{3p_{1}-p_{2}}{2p_{1}p_{2}}\right)  \int\limits_{\pi/2}^{{\pi}}%
\frac{\sqrt{2}q\left(  \tau\right)  }{2}\cos\left(  \frac{s(2\tau-\Delta
(\tau))}{p_{2}}+\frac{\pi}{4}\right)  d\tau\right.
\end{align*}%
\begin{align}
&  -\cos\left(  s\pi\frac{3p_{1}-p_{2}}{2p_{1}p_{2}}\right)  \int
\limits_{\pi/2}^{{\pi}}\frac{\sqrt{2}q\left(  \tau\right)  }{2}\sin\left(
\frac{s(2\tau-\Delta(\tau))}{p_{2}}+\frac{\pi}{4}\right)  d\tau\nonumber\\
&  +\sin\left(  \frac{s\pi\left(  p_{1}+p_{2}\right)  }{2p_{1}p_{2}}\right)
\int\limits_{\pi/2}^{{\pi}}\frac{\sqrt{2}q\left(  \tau\right)  }{2}\cos\left(
\frac{s\Delta(\tau)}{p_{2}}-\frac{\pi}{4}\right)  d\tau\nonumber\\
&  \left.  -\cos\left(  \frac{s\pi\left(  p_{1}+p_{2}\right)  }{2p_{1}p_{2}%
}\right)  \int\limits_{\pi/2}^{{\pi}}\frac{\sqrt{2}q\left(  \tau\right)  }%
{2}\sin\left(  \frac{s\Delta(\tau)}{p_{2}}-\frac{\pi}{4}\right)  d\tau\right]
+O\left(  \frac{1}{s}\right)  \tag{31}%
\end{align}

\noindent Let%
\begin{align}
A\left(  x,\>s,\>\Delta(\tau)\right)   &  =\int\limits_{0}^{x}\frac{\sqrt
{2}q\left(  \tau\right)  }{2}\sin\left(  \frac{s\Delta(\tau)}{p_{1}}-\frac
{\pi}{4}\right)  d\tau,\nonumber\\
B(x,\>s,\>\Delta(\tau))  &  =\int\limits_{0}^{x}\frac{\sqrt{2}q\left(
\tau\right)  }{2}\cos\left(  \frac{s\Delta(\tau)}{p_{1}}-\frac{\pi}{4}\right)
d\tau) \tag{32}%
\end{align}
\noindent It is obviously that these functions are bounded for $\>0\leq
x\leq\frac{\pi}{2},\>0<s<\infty$. Let%
\begin{align}
C\left(  x,\>s,\>\Delta(\tau)\right)   &  =\int\limits_{\pi/2}^{{x}}%
\frac{\sqrt{2}q\left(  \tau\right)  }{2}\sin\left(  \frac{s\Delta(\tau)}%
{p_{2}}-\frac{\pi}{4}\right)  d\tau,\nonumber\\
D\left(  x,\>s,\>\Delta(\tau)\right)   &  =\int\limits_{\pi/2}^{{x}}%
\frac{\sqrt{2}q\left(  \tau\right)  }{2}\cos\left(  \frac{s\Delta(\tau)}%
{p_{2}}-\frac{\pi}{4}\right)  d\tau\tag{33}%
\end{align}

\noindent It is obviously that these functions are bounded for $\>\frac{\pi
}{2}\leq x\leq\pi,\>0<s<\infty$.

\noindent Under the conditions a) and b) the following formulas%
\begin{align}
\int\limits_{0}^{x}\frac{\sqrt{2}q\left(  \tau\right)  }{2}\cos\left(
\frac{s\left(  2\tau-\Delta(\tau)\right)  }{p_{1}}+\frac{\pi}{4}\right)  d\tau
&  =O\left(  \frac{1}{s}\right)  ,\nonumber\\
\int\limits_{0}^{x}\frac{\sqrt{2}q\left(  \tau\right)  }{2}\sin\left(
\frac{s\left(  2\tau-\Delta(\tau)\right)  }{p_{1}}+\frac{\pi}{4}\right)  d\tau
&  =O\left(  \frac{1}{s}\right)  ,\nonumber\\
\int\limits_{\pi/2}^{{x}}\frac{\sqrt{2}q\left(  \tau\right)  }{2}\cos\left(
\frac{s(2\tau-\Delta(\tau))}{p_{2}}+\frac{\pi}{4}\right)  d\tau &  =O\left(
\frac{1}{s}\right)  ,\nonumber\\
\int\limits_{\pi/2}^{{x}}\frac{\sqrt{2}q\left(  \tau\right)  }{2}\sin\left(
\frac{s(2\tau-\Delta(\tau))}{p_{2}}+\frac{\pi}{4}\right)  d\tau &  =O\left(
\frac{1}{s}\right)  . \tag{34}%
\end{align}
\noindent can be proved by the same technique in Lemma 3.3.3 in \cite{Nk}.
From $(31),(32),(33),(34)$ we have%
\[
\cot\left(  \frac{s\pi\left(  p_{1}+p_{2}\right)  }{2p_{1}p_{2}}+\frac{\pi}%
{4}\right)  =\frac{1}{s}\left[  \frac{\gamma_{2}}{\delta_{2}}+\frac
{d\gamma_{1}B\left(  x,s,\Delta\left(  \tau\right)  \right)  }{p_{1}\delta
_{1}}+\frac{d\gamma_{1}D\left(  x,s,\Delta\left(  \tau\right)  \right)
}{p_{2}\delta_{1}}\right]  +O\left(  \frac{1}{s^{2}}\right)
\]
Now using $s_{n}=\frac{p_{1}p_{2}\left(  4n-3\right)  }{2\left(  p_{1}%
+p_{2}\right)  }+\delta_{n}\>$we get%
\[
\cot\left(  \frac{\left(  2n-1\right)  \pi}{2}+\frac{\pi\left(  p_{1}%
+p_{2}\right)  \delta_{n}}{2p_{1}p_{2}}\right)  =-\tan\frac{\pi\left(
p_{1}+p_{2}\right)  \delta_{n}}{2p_{1}p_{2}}=\frac{2\left(  p_{1}%
+p_{2}\right)  }{p_{1}p_{2}\left(  4n-3\right)  }%
\]%
\[
\times\left[  \frac{\gamma_{2}}{\delta_{2}}+\frac{d\gamma_{1}B\left(
\frac{\pi}{2},\frac{p_{1}p_{2}\left(  4n-3\right)  }{2\left(  p_{1}%
+p_{2}\right)  },\Delta\left(  \tau\right)  \right)  }{p_{1}\delta_{1}}%
+\frac{d\gamma_{1}D\left(  \pi,\frac{p_{1}p_{2}\left(  4n-3\right)  }{2\left(
p_{1}+p_{2}\right)  },\Delta\left(  \tau\right)  \right)  }{p_{2}\delta_{1}%
}\right]  +O\left(  \frac{1}{n^{2}}\right)
\]

\noindent and finally%
\begin{align*}
\delta_{n}  &  =\frac{4}{\left(  4n-3\right)  \pi}\left[  \frac{\gamma_{2}%
}{\delta_{2}}+\frac{d\gamma_{1}B\left(  \frac{\pi}{2},\frac{p_{1}p_{2}\left(
4n-3\right)  }{2\left(  p_{1}+p_{2}\right)  },\Delta\left(  \tau\right)
\right)  }{p_{1}\delta_{1}}+\frac{d\gamma_{1}D\left(  \pi,\frac{p_{1}%
p_{2}\left(  4n-3\right)  }{2\left(  p_{1}+p_{2}\right)  },\Delta\left(
\tau\right)  \right)  }{p_{2}\delta_{1}}\right] \\
&  +O\left(  \frac{1}{n^{2}}\right)
\end{align*}
Thus%
\[
s_{n}=\frac{p_{1}p_{2}\left(  4n-3\right)  }{2\left(  p_{1}+p_{2}\right)
}+\frac{4}{\left(  4n-3\right)  \pi}\left[  \frac{\gamma_{2}}{\delta_{2}%
}+\frac{d\gamma_{1}B\left(  \frac{\pi}{2},\frac{p_{1}p_{2}\left(  4n-3\right)
}{2\left(  p_{1}+p_{2}\right)  },\Delta\left(  \tau\right)  \right)  }%
{p_{1}\delta_{1}}\right.
\]%
\begin{equation}
\left.  +\frac{d\gamma_{1}D\left(  \pi,\frac{p_{1}p_{2}\left(  4n-3\right)
}{2\left(  p_{1}+p_{2}\right)  },\Delta\left(  \tau\right)  \right)  }%
{p_{2}\delta_{1}}\right]  +O\left(  \frac{1}{n^{2}}\right)  \tag{35}%
\end{equation}
Thus, we have proven the following theorem.

\textbf{Theorem 4} \ If conditions a) and b) are satisfied then, the positive
eigenvalues $\lambda_{n}=s_{n}^{2}\>$ of the problem (1)-(5) have the $(35)$
asymptotic representation for$\>n\rightarrow\infty\>$.

We now may obtain a sharper asymptotic formula for the eigenfunctions. From
$(8)$ and $(29)$%
\begin{align*}
w_{1}(x,\lambda)  &  =\sqrt{2}p_{1}\cos\left(  \frac{sx}{p_{1}}+\frac{\pi}%
{4}\right) \\
&  -\frac{\sqrt{2}}{sp_{1}}\int\limits_{0}^{{x}}q\left(  \tau\right)
\sin\left(  \frac{s}{p_{1}}\left(  x-\tau\right)  \right)  \cos\left(
\frac{s\left(  \tau-\Delta\left(  \tau\right)  \right)  }{p_{1}}+\frac{\pi}%
{4}\right)  d\tau+O\left(  \frac{1}{s^{2}}\right)  ,
\end{align*}%
\begin{align*}
w_{1}(x,\lambda)  &  =\sqrt{2}p_{1}\cos\left(  \frac{sx}{p_{1}}+\frac{\pi}%
{4}\right)  -\frac{\sqrt{2}}{2sp_{1}}\int\limits_{0}^{{x}}q\left(
\tau\right)  \left[  \left[  \sin\left(  \frac{sx}{p_{1}}-\frac{\pi}%
{4}\right)  \cos\left(  \frac{s\left(  2\tau-\Delta\left(  \tau\right)
\right)  }{p_{1}}\right)  \right.  \right. \\
&  \left.  -\cos\left(  \frac{sx}{p_{1}}-\frac{\pi}{4}\right)  \sin\left(
\frac{s\left(  2\tau-\Delta\left(  \tau\right)  \right)  }{p_{1}}\right)
\right]  +\left[  \sin\left(  \frac{sx}{p_{1}}+\frac{\pi}{4}\right)
\cos\left(  \frac{s\Delta\left(  \tau\right)  }{p_{1}}\right)  \right. \\
&  \left.  \left.  -\cos\left(  \frac{sx}{p_{1}}+\frac{\pi}{4}\right)
\sin\left(  \frac{s\Delta\left(  \tau\right)  }{p_{1}}\right)  \right]
\right]  +O\left(  \frac{1}{s^{2}}\right)  .
\end{align*}
Thus, from $(32),(33)$ and $(34)$%
\[
w_{1}(x,\lambda)=\cos\left(  \frac{sx}{p_{1}}+\frac{\pi}{4}\right)  \left[
\sqrt{2}p_{1}+\frac{A\left(  x,\>s,\>\Delta(\tau)\right)  }{sp_{1}}\right]
\]%
\begin{equation}
-\frac{\sin\left(  \frac{sx}{p_{1}}+\frac{\pi}{4}\right)  B\left(
x,\>s,\>\Delta(\tau)\right)  }{sp_{1}}+O\left(  \frac{1}{s^{2}}\right)  .
\tag{36}%
\end{equation}
\noindent Replacing$\>s$ by$\>s_{n}$ and using $(35)$ for $x\in\left[
0,\frac{\pi}{2}\right)  $ we have
\begin{align}
u_{1n}(x)  &  =\cos\left(  \frac{p_{2}\left(  4n-3\right)  x}{2\left(
p_{1}+p_{2}\right)  }+\frac{\pi}{4}\right)  \left[  \sqrt{2}p_{1}%
+\frac{2\left(  p_{1}+p_{2}\right)  A\left(  x,\>\frac{p_{1}p_{2}\left(
4n-3\right)  }{2\left(  p_{1}+p_{2}\right)  },\>\Delta(\tau)\right)  }%
{p_{1}^{2}p_{2}(4n-3)}\right] \nonumber\\
&  -\sin\left(  \frac{p_{2}\left(  4n-3\right)  x}{2\left(  p_{1}%
+p_{2}\right)  }+\frac{\pi}{4}\right)  \left[  \frac{4\sqrt{2}}{(4n-3)\pi
}\left[  \frac{p_{1}\gamma_{2}}{\delta_{2}}+\frac{d\gamma_{1}B\left(
\frac{\pi}{2},\frac{p_{1}p_{2}\left(  4n-3\right)  }{2\left(  p_{1}%
+p_{2}\right)  },\Delta\left(  \tau\right)  \right)  }{\delta_{1}}\right.
\right. \nonumber\\
&  \left.  \left.  +\frac{dp_{1}\gamma_{1}D\left(  \pi,\frac{p_{1}p_{2}\left(
4n-3\right)  }{2\left(  p_{1}+p_{2}\right)  },\Delta\left(  \tau\right)
\right)  }{p_{2}\delta_{1}}\right]  \right]  +O\left(  \frac{1}{n^{2}}\right)
. \tag{37}%
\end{align}
From $(16),(29)$ and $(32)$, we have%
\begin{align}
\frac{w_{1}^{^{\prime}}\left(  x,\lambda\right)  }{s}  &  =-\sqrt{2}%
\sin\left(  \frac{sx}{p_{1}}+\frac{\pi}{4}\right)  -\frac{B\left(
x,\>s,\>\Delta(\tau)\right)  \cos\left(  \frac{sx}{p_{1}}+\frac{\pi}%
{4}\right)  }{p_{1}^{2}s}\nonumber\\
&  -\frac{A\left(  x,\>s,\>\Delta(\tau)\right)  \sin\left(  \frac{sx}{p_{1}%
}+\frac{\pi}{4}\right)  }{p_{1}^{2}s}+O\left(  \frac{1}{s^{2}}\right)  \text{
},x\in\left(  0,\frac{\pi}{2}\right]  \tag{38}%
\end{align}
From $(9),(30),(34),(36)$ and $(38)$ we have%
\begin{align*}
w_{2}\left(  x,\lambda\right)   &  =\frac{\gamma_{1}}{2\delta_{1}}\left\{
\left[  \sqrt{2}p_{1}+\frac{A\left(  \frac{\pi}{2},s,\Delta(\tau)\right)
}{sp_{1}}\right]  \right. \\
&  \times\left[  \cos\left(  s\left(  \frac{\pi(p_{1}+p_{2})}{2p_{1}p_{2}%
}-\frac{x}{p_{2}}\right)  +\frac{\pi}{4}\right)  +\cos\left(  s\left(
\frac{\pi(p_{2}-p_{1})}{2p_{1}p_{2}}+\frac{x}{p_{2}}\right)  +\frac{\pi}%
{4}\right)  \right] \\
&  -\left[  \sin\left(  s\left(  \frac{\pi(p_{1}+p_{2})}{2p_{1}p_{2}}-\frac
{x}{p_{2}}\right)  +\frac{\pi}{4}\right)  +\sin\left(  s\left(  \frac
{\pi(p_{2}-p_{1})}{2p_{1}p_{2}}+\frac{x}{p_{2}}\right)  +\frac{\pi}{4}\right)
\right] \\
&  \times\left.  \frac{B\left(  \frac{\pi}{2},s,\Delta(\tau)\right)  }{sp_{1}%
}\right\}  -\left\{  \frac{\gamma_{2}p_{2}\sqrt{2}}{2\delta_{2}}+\frac
{\gamma_{2}p_{2}}{2\delta_{2}p_{1}^{2}s}\right\} \\
&  \times\left\{  -\cos\left(  s\left(  \frac{\pi(p_{2}-p_{1})}{2p_{1}p_{2}%
}+\frac{x}{p_{2}}\right)  +\frac{\pi}{4}\right)  +\cos\left(  s\left(
\frac{\pi(p_{1}+p_{2})}{2p_{1}p_{2}}-\frac{x}{p_{2}}\right)  +\frac{\pi}%
{4}\right)  \right\} \\
&  -\frac{\gamma_{2}p_{2}B\left(  \frac{\pi}{2},s,\Delta(\tau)\right)
}{2\delta_{2}p_{1}^{2}s}\\
&  \times\left\{  \sin\left(  s\left(  \frac{\pi(p_{2}-p_{1})}{2p_{1}p_{2}%
}+\frac{x}{p_{2}}\right)  +\frac{\pi}{4}\right)  -\sin\left(  s\left(
\frac{\pi(p_{1}+p_{2})}{2p_{1}p_{2}}-\frac{x}{p_{2}}\right)  +\frac{\pi}%
{4}\right)  \right\} \\
&  -\frac{\gamma_{1}p_{1}}{sp_{2}\delta_{1}}\left\{  D\left(  x,s,\Delta
(\tau)\right)  \sin s\left(  \frac{x}{p_{2}}+\frac{\pi(p_{2}-p_{1})}%
{2p_{1}p_{2}}\right)  \right. \\
&  \left.  -C\left(  x,s,\Delta(\tau)\right)  \cos s\left(  \frac{x}{p_{2}%
}+\frac{\pi(p_{2}-p_{1})}{2p_{1}p_{2}}\right)  \right\}  +O\left(  \frac
{1}{s^{2}}\right)  ,\text{ \ \ }x\in\left(  \frac{\pi}{2},\pi\right]  .
\end{align*}
\noindent Now, replacing $\>s\>$ by $\>s_{n}\>$and using $(35)$, we have
\[
u_{2n}(x)=\frac{\gamma_{1}}{2\delta_{1}}\left\{  \left[  (-1)^{n+1}\sin\left(
\frac{(4n-3)p_{1}x}{2\left(  p_{1}+p_{2}\right)  }\right)  +\cos\left(
\frac{(4n-3)\left(  p_{2}-p_{1}\right)  \pi}{4\left(  p_{1}+p_{2}\right)
}+\frac{(4n-3)p_{1}x}{2\left(  p_{1}+p_{2}\right)  }+\frac{\pi}{4}\right)
\right]  \right.
\]%
\begin{align}
&  \,\times\left[  \sqrt{2}p_{1}+\frac{2\left(  p_{1}+p_{2}\right)  A\left(
\frac{\pi}{2},\frac{(4n-3)p_{1}p_{2}}{2\left(  p_{1}+p_{2}\right)  }%
,\Delta(\tau))\right)  }{(4n-3)p_{1}^{2}p_{2}}\right]  +\left[  (-1)^{n}%
\cos\left(  \frac{(4n-3)p_{1}x}{2\left(  p_{1}+p_{2}\right)  }\right)
-\right. \\
&  \left.  \left.  \sin\left(  \frac{(4n-3)\left(  p_{2}-p_{1}\right)  \pi
}{4\left(  p_{1}+p_{2}\right)  }+\frac{(4n-3)p_{1}x}{2\left(  p_{1}%
+p_{2}\right)  }+\frac{\pi}{4}\right)  \right]  \left[  \frac{2\left(
p_{1}+p_{2}\right)  B\left(  \frac{\pi}{2},\frac{(4n-3)p_{1}p_{2}}{2\left(
p_{1}+p_{2}\right)  },\Delta(\tau))\right)  }{(4n-3)p_{1}^{2}p_{2}}\right]
\right\} \nonumber\\
&  -\frac{\sqrt{2}\gamma_{2}p_{2}}{2\delta_{2}}\left\{  (-1)^{n+1}\sin\left(
\frac{(4n-3)p_{1}x}{2\left(  p_{1}+p_{2}\right)  }\right)  +\frac{4}%
{(4n-3)\pi}\left[  \,(-1)^{n}\cos\left(  \frac{(4n-3)p_{1}x}{2\left(
p_{1}+p_{2}\right)  }\right)  +\right.  \right. \nonumber\\
&  \left.  \times\sin\left(  \frac{(4n-3)\left(  p_{2}-p_{1}\right)  \pi
}{4\left(  p_{1}+p_{2}\right)  }+\frac{(4n-3)p_{1}x}{2\left(  p_{1}%
+p_{2}\right)  }+\frac{\pi}{4}\right)  \right]  \,\left[  \frac{\gamma_{2}%
}{\delta_{2}}+\frac{d\gamma_{1}B\left(  \frac{\pi}{2},\frac{(4n-3)p_{1}p_{2}%
}{2\left(  p_{1}+p_{2}\right)  },\Delta(\tau))\right)  }{p_{1}\delta_{1}%
}+\right. \nonumber\\
&  \,\left.  \left.  \frac{d\gamma_{1}D\left(  \pi,\frac{(4n-3)p_{1}p_{2}%
}{2\left(  p_{1}+p_{2}\right)  },\Delta(\tau))\right)  }{p_{2}\delta_{1}%
}\right]  +\cos\left(  \frac{(4n-3)\left(  p_{2}-p_{1}\right)  \pi}{4\left(
p_{1}+p_{2}\right)  }+\frac{(4n-3)p_{1}x}{2\left(  p_{1}+p_{2}\right)  }%
+\frac{\pi}{4}\right)  \right\} \nonumber\\
&  +\frac{\gamma_{2}\left(  p_{1}+p_{2}\right)  }{\delta_{2}p_{1}^{3}%
(4n-3)}B\left(  \frac{\pi}{2},\frac{(4n-3)p_{1}p_{2}}{2\left(  p_{1}%
+p_{2}\right)  },\Delta(\tau))\right)  \left\{  (-1)^{n}\cos\left(
\frac{(4n-3)p_{1}x}{2\left(  p_{1}+p_{2}\right)  }\right)  \right. \nonumber\\
&  \left.  -\sin\left(  \frac{(4n-3)\left(  p_{2}-p_{1}\right)  \pi}{4\left(
p_{1}+p_{2}\right)  }+\frac{(4n-3)p_{1}x}{2\left(  p_{1}+p_{2}\right)  }%
+\frac{\pi}{4}\right)  \right\}  +\frac{2\gamma_{2}\left(  p_{1}+p_{2}\right)
A\left(  \frac{\pi}{2},\frac{(4n-3)p_{1}p_{2}}{2\left(  p_{1}+p_{2}\right)
},\Delta(\tau)\right)  }{\delta_{2}p_{1}^{3}(4n-3)}\nonumber\\
&  \times\left\{  (-1)^{n}\sin\left(  \frac{(4n-3)p_{1}x}{2\left(  p_{1}%
+p_{2}\right)  }\right)  +\cos\left(  \frac{(4n-3)\left(  p_{2}-p_{1}\right)
\pi}{4\left(  p_{1}+p_{2}\right)  }+\frac{(4n-3)p_{1}x}{2\left(  p_{1}%
+p_{2}\right)  }+\frac{\pi}{4}\right)  \right\}  -\nonumber\\
&  \frac{2\gamma_{1}\left(  p_{1}+p_{2}\right)  }{(4n-3)\delta_{1}p_{2}^{2}%
}\left\{  \sin\left(  \frac{(4n-3)\left(  p_{2}-p_{1}\right)  \pi}{4\left(
p_{1}+p_{2}\right)  }+\frac{(4n-3)p_{1}x}{2\left(  p_{1}+p_{2}\right)  }%
+\frac{\pi}{4}\right)  D\left(  x,\frac{(4n-3)p_{1}p_{2}}{2\left(  p_{1}%
+p_{2}\right)  },\Delta(\tau)\right)  \right. \nonumber\\
&  \left.  -C\left(  x,\frac{(4n-3)p_{1}p_{2}}{2\left(  p_{1}+p_{2}\right)
},\Delta(\tau)\right)  \cos\left(  \frac{(4n-3)\left(  p_{2}-p_{1}\right)
\pi}{4\left(  p_{1}+p_{2}\right)  }+\frac{(4n-3)p_{1}x}{2\left(  p_{1}%
+p_{2}\right)  }+\frac{\pi}{4}\right)  \right\}  +O\bigl({\frac{1}{{n^{2}}}%
}\bigr ). \tag{39}%
\end{align}
Thus, we have proven the following theorem.

\textbf{Theorem 5} \ If conditions a) and b) are satisfied then, the
eigenfunctions$\>u_{n}(x)\>$of the problem (1)-(5) have the following
asymptotic representation for$\>n\rightarrow\infty\>$:
\[
u_{n}(x)=\left\{
\begin{array}
[c]{lll}%
u_{1n}(x) & \mbox{for} & x\in\lbrack0,{\frac{\pi}{2}})\\
&  & \\
u_{2n}(x) & \mbox{for} & x\in({\frac{\pi}{2}},\pi]
\end{array}
\right.
\]
\noindent where$\>u_{1n}(x)\>$and$\>u_{2n}(x)\>$defined as in (37) and (39) respectively.

\end{document}